\documentclass[a4paper,twoside,11pt,english]{article}

\usepackage[all]{xy}
\usepackage[intlimits]{amsmath}
\usepackage{amssymb,amsthm}
\usepackage{graphicx}
\usepackage[dvips,dvipdf,pdftex]{epsfig}
\usepackage{babel}
\usepackage[T1]{fontenc}
\usepackage{fancyhdr,stmaryrd,color}
\usepackage{times,euler,euscript,eufrak}

\DeclareGraphicsExtensions{.eps}

  \newtheorem{thm}{Theorem}[section]

  \newtheorem{prop}[thm]{Proposition}
\theoremstyle{definition}
  \newtheorem{defn}[thm]{Definition}
  
\theoremstyle{remark}

\DeclareMathAlphabet{\mathscr}{T1}{pzc}{m}{it}
\newcommand{\Nb}{\mathbb{N}}
\newcommand{\lz}{\mathrm{L}(\mathbb{Z}_2)}
\newcommand{\zb}{\mathbb{Z}/2\mathbb{Z}}
\renewcommand{\phi}{\varphi}
\renewcommand{\epsilon}{\varepsilon}
\newcommand{\ul}[1]{\underline{#1}}
\newcommand{\fl}{\rightarrow}
\newcommand{\lsur}{\longmapsto}
\newcommand{\red}[1]{\rightarrow\!\!_{{\scriptscriptstyle #1}}}
\newcommand{\mon}[1]{\langle #1 \rangle}

\DeclareMathOperator{\Id}{Id}

\newcommand{\emptysectionmark}[1]{\markboth{\textbf{#1}}{\textbf{#1}}}

\fancypagestyle{plain}{
  \fancyhf{}\fancyfoot[LC,RC]{}
  \fancyfoot[LE,RO]{\textbf{\thepage}}
  
  }

\begin{document}
\thispagestyle{plain}
\hfill {\large \textbf{4th May 2005 -- Modified: 18th November 2005}}

\vspace{0.5mm}
\hrule height 1.5pt
\vspace{3mm}
{\LARGE \textbf{TERMINATION ORDERS FOR $\mathbf{3}$-POLYGRAPHS}} 

\vspace{1.5mm}
\indent{\LARGE \textbf{Yves Guiraud\footnote{Institut de mathématiques de Luminy, Marseille, France -- guiraud@iml.univ-mrs.fr}}}
\vspace{2.5mm}
\hrule height 1.5pt

\bigskip
\begin{minipage}{140mm}
\textbf{Résumé :} Cette note présente la première classe connue d'ordres de terminaison adaptés aux $3$-polygraphes, ainsi qu'une application. 

\smallskip
\noindent \textbf{Abstract:} This note presents the first known class of termination orders for $3$-polygraphs, together with an application.
\end{minipage}

\bigskip
\noindent Polygraphs are cellular presentations of higher-dimensional categories introduced in [Burroni 1993]. They have been proved to generalize term rewriting systems but they lack some tools widely used in the field. This note presents a result developped in [Guiraud 2004] which fills this gap for some $3$-dimensional polygraphs: it introduces a method to craft \emph{termination orders}, one of the most useful ways to prove that computations specified by a formal system always end after a finite number of transformations. 

\section{Notions about $\mathbf{3}$-polygraphs}

\noindent The formal definition of polygraphs can be found in [Burroni 1993]. Here, we restrict ourselves to the case of a \emph{$\mathit{2}$-polygraph with one $\mathit{0}$-cell and one $\mathit{1}$-cell}: this is a graph $\Sigma$ over the set of natural numbers. Elements of $\Sigma$ are called \emph{$2$-dimensional cells} or \emph{circuits}. Two $2$-cells are \emph{parallel} when thay have the same source and the same target. A $2$-dimensional cell $\phi:m\fl n$ is graphically pictured as a circuit with $m$ inputs and $n$ outputs:
\begin{center}
\begin{picture}(0,0)%
\includegraphics{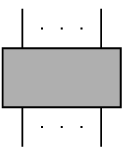}%
\end{picture}%
\setlength{\unitlength}{4144sp}%
\begingroup\makeatletter\ifx\SetFigFont\undefined%
\gdef\SetFigFont#1#2#3#4#5{%
  \reset@font\fontsize{#1}{#2pt}%
  \fontfamily{#3}\fontseries{#4}\fontshape{#5}%
  \selectfont}%
\fi\endgroup%
\begin{picture}(564,958)(79,-56)
\put(361,-16){\makebox(0,0)[b]{\smash{{\SetFigFont{10}{12.0}{\rmdefault}{\mddefault}{\updefault}{\color[rgb]{0,0,0}$n$}%
}}}}
\put(361,389){\makebox(0,0)[b]{\smash{{\SetFigFont{10}{12.0}{\rmdefault}{\mddefault}{\updefault}{\color[rgb]{0,0,0}$\phi$}%
}}}}
\put(361,794){\makebox(0,0)[b]{\smash{{\SetFigFont{10}{12.0}{\rmdefault}{\mddefault}{\updefault}{\color[rgb]{0,0,0}$m$}%
}}}}
\end{picture}%
\end{center}

\noindent Given such a $2$-polygraph $\Sigma$, one builds another $2$-polygraph $\mon{\Sigma}$: its $2$-cells are all the circuits one can build from the ones in $\Sigma$, by either (horizontal) juxtaposition or (vertical) plugging. These two operations are pictured this way:
\begin{center}
\begin{picture}(0,0)%
\includegraphics{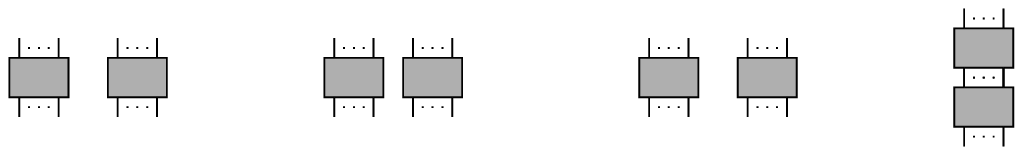}%
\end{picture}%
\setlength{\unitlength}{4144sp}%
\begingroup\makeatletter\ifx\SetFigFont\undefined%
\gdef\SetFigFont#1#2#3#4#5{%
  \reset@font\fontsize{#1}{#2pt}%
  \fontfamily{#3}\fontseries{#4}\fontshape{#5}%
  \selectfont}%
\fi\endgroup%
\begin{picture}(4737,654)(-44,-568)
\put(1261,-151){\makebox(0,0)[b]{\smash{{\SetFigFont{10}{12.0}{\rmdefault}{\mddefault}{\updefault}{\color[rgb]{0,0,0}$\star_0$}%
}}}}
\put(-44,-286){\makebox(0,0)[b]{\smash{{\SetFigFont{10}{12.0}{\rmdefault}{\mddefault}{\updefault}{\color[rgb]{0,0,0}$\Big($}%
}}}}
\put(946,-286){\makebox(0,0)[b]{\smash{{\SetFigFont{10}{12.0}{\rmdefault}{\mddefault}{\updefault}{\color[rgb]{0,0,0}$\Big)$}%
}}}}
\put(451,-286){\makebox(0,0)[b]{\smash{{\SetFigFont{10}{12.0}{\rmdefault}{\mddefault}{\updefault}{\color[rgb]{0,0,0}$,$}%
}}}}
\put(2836,-286){\makebox(0,0)[b]{\smash{{\SetFigFont{10}{12.0}{\rmdefault}{\mddefault}{\updefault}{\color[rgb]{0,0,0}$\Big($}%
}}}}
\put(3826,-286){\makebox(0,0)[b]{\smash{{\SetFigFont{10}{12.0}{\rmdefault}{\mddefault}{\updefault}{\color[rgb]{0,0,0}$\Big)$}%
}}}}
\put(3331,-286){\makebox(0,0)[b]{\smash{{\SetFigFont{10}{12.0}{\rmdefault}{\mddefault}{\updefault}{\color[rgb]{0,0,0}$,$}%
}}}}
\put(226,-286){\makebox(0,0)[b]{\smash{{\SetFigFont{10}{12.0}{\rmdefault}{\mddefault}{\updefault}{\color[rgb]{0,0,0}$f$}%
}}}}
\put(676,-286){\makebox(0,0)[b]{\smash{{\SetFigFont{10}{12.0}{\rmdefault}{\mddefault}{\updefault}{\color[rgb]{0,0,0}$g$}%
}}}}
\put(1666,-286){\makebox(0,0)[b]{\smash{{\SetFigFont{10}{12.0}{\rmdefault}{\mddefault}{\updefault}{\color[rgb]{0,0,0}$f$}%
}}}}
\put(2026,-286){\makebox(0,0)[b]{\smash{{\SetFigFont{10}{12.0}{\rmdefault}{\mddefault}{\updefault}{\color[rgb]{0,0,0}$g$}%
}}}}
\put(3106,-286){\makebox(0,0)[b]{\smash{{\SetFigFont{10}{12.0}{\rmdefault}{\mddefault}{\updefault}{\color[rgb]{0,0,0}$f$}%
}}}}
\put(4546,-151){\makebox(0,0)[b]{\smash{{\SetFigFont{10}{12.0}{\rmdefault}{\mddefault}{\updefault}{\color[rgb]{0,0,0}$f$}%
}}}}
\put(3556,-286){\makebox(0,0)[b]{\smash{{\SetFigFont{10}{12.0}{\rmdefault}{\mddefault}{\updefault}{\color[rgb]{0,0,0}$g$}%
}}}}
\put(4546,-421){\makebox(0,0)[b]{\smash{{\SetFigFont{10}{12.0}{\rmdefault}{\mddefault}{\updefault}{\color[rgb]{0,0,0}$g$}%
}}}}
\put(1261,-286){\makebox(0,0)[b]{\smash{{\SetFigFont{10}{12.0}{\rmdefault}{\mddefault}{\updefault}{\color[rgb]{0,0,0}$\lsur$}%
}}}}
\put(4141,-286){\makebox(0,0)[b]{\smash{{\SetFigFont{10}{12.0}{\rmdefault}{\mddefault}{\updefault}{\color[rgb]{0,0,0}$\lsur$}%
}}}}
\put(4141,-151){\makebox(0,0)[b]{\smash{{\SetFigFont{10}{12.0}{\rmdefault}{\mddefault}{\updefault}{\color[rgb]{0,0,0}$\star_1$}%
}}}}
\end{picture}%
\end{center}

\noindent These constructions are considered \emph{modulo isotopy} (or homeomorphic deformation):
\begin{center}
\begin{picture}(0,0)%
\includegraphics{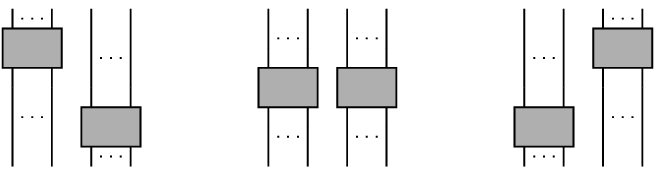}%
\end{picture}%
\setlength{\unitlength}{4144sp}%
\begingroup\makeatletter\ifx\SetFigFont\undefined%
\gdef\SetFigFont#1#2#3#4#5{%
  \reset@font\fontsize{#1}{#2pt}%
  \fontfamily{#3}\fontseries{#4}\fontshape{#5}%
  \selectfont}%
\fi\endgroup%
\begin{picture}(2994,744)(79,17)
\put(2161,344){\makebox(0,0)[b]{\smash{{\SetFigFont{10}{12.0}{\rmdefault}{\mddefault}{\updefault}{\color[rgb]{0,0,0}$\equiv$}%
}}}}
\put(226,524){\makebox(0,0)[b]{\smash{{\SetFigFont{10}{12.0}{\rmdefault}{\mddefault}{\updefault}{\color[rgb]{0,0,0}$f$}%
}}}}
\put(586,164){\makebox(0,0)[b]{\smash{{\SetFigFont{10}{12.0}{\rmdefault}{\mddefault}{\updefault}{\color[rgb]{0,0,0}$g$}%
}}}}
\put(1396,344){\makebox(0,0)[b]{\smash{{\SetFigFont{10}{12.0}{\rmdefault}{\mddefault}{\updefault}{\color[rgb]{0,0,0}$f$}%
}}}}
\put(1756,344){\makebox(0,0)[b]{\smash{{\SetFigFont{10}{12.0}{\rmdefault}{\mddefault}{\updefault}{\color[rgb]{0,0,0}$g$}%
}}}}
\put(2566,164){\makebox(0,0)[b]{\smash{{\SetFigFont{10}{12.0}{\rmdefault}{\mddefault}{\updefault}{\color[rgb]{0,0,0}$f$}%
}}}}
\put(2926,524){\makebox(0,0)[b]{\smash{{\SetFigFont{10}{12.0}{\rmdefault}{\mddefault}{\updefault}{\color[rgb]{0,0,0}$g$}%
}}}}
\put(991,344){\makebox(0,0)[b]{\smash{{\SetFigFont{10}{12.0}{\rmdefault}{\mddefault}{\updefault}{\color[rgb]{0,0,0}$\equiv$}%
}}}}
\end{picture}%
\end{center}

\begin{defn}
A \emph{$\mathit{3}$-polygraph with one $\mathit{0}$-cell and one $\mathit{1}$-cell} is a pair $(\Sigma,R)$ such that $\Sigma$ is a $2$-polygraph with one $0$-cell and one $1$-cell and $R$ is a graph over $\mon{\Sigma}$ made of arrows between parallel circuits. An element of $R$ is called a \emph{$\mathit{3}$-cell}.

The \emph{reduction relation generated by $R$} is the binary relation on circuits of $\mon{\Sigma}$ defined by $f\red{R}\:g$ whenever there exists a $3$-cell $\alpha:f_0\fl g_0$, together with two circuits $h$ and $k$, such that the following relations have a meaning and hold:
\begin{center}
\begin{picture}(0,0)%
\includegraphics{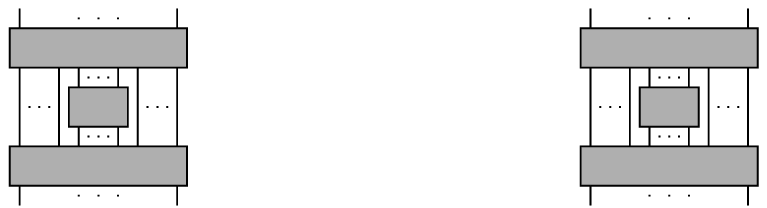}%
\end{picture}%
\setlength{\unitlength}{4144sp}%
\begingroup\makeatletter\ifx\SetFigFont\undefined%
\gdef\SetFigFont#1#2#3#4#5{%
  \reset@font\fontsize{#1}{#2pt}%
  \fontfamily{#3}\fontseries{#4}\fontshape{#5}%
  \selectfont}%
\fi\endgroup%
\begin{picture}(4035,924)(388,-163)
\put(496,254){\makebox(0,0)[b]{\smash{{\SetFigFont{10}{12.0}{\rmdefault}{\mddefault}{\updefault}{\color[rgb]{0,0,0}$f$}%
}}}}
\put(1396,524){\makebox(0,0)[b]{\smash{{\SetFigFont{10}{12.0}{\rmdefault}{\mddefault}{\updefault}{\color[rgb]{0,0,0}$h$}%
}}}}
\put(1396,-16){\makebox(0,0)[b]{\smash{{\SetFigFont{10}{12.0}{\rmdefault}{\mddefault}{\updefault}{\color[rgb]{0,0,0}$k$}%
}}}}
\put(4006,524){\makebox(0,0)[b]{\smash{{\SetFigFont{10}{12.0}{\rmdefault}{\mddefault}{\updefault}{\color[rgb]{0,0,0}$h$}%
}}}}
\put(4006,-16){\makebox(0,0)[b]{\smash{{\SetFigFont{10}{12.0}{\rmdefault}{\mddefault}{\updefault}{\color[rgb]{0,0,0}$k$}%
}}}}
\put(1396,254){\makebox(0,0)[b]{\smash{{\SetFigFont{10}{12.0}{\rmdefault}{\mddefault}{\updefault}{\color[rgb]{0,0,0}$f_0$}%
}}}}
\put(4006,254){\makebox(0,0)[b]{\smash{{\SetFigFont{10}{12.0}{\rmdefault}{\mddefault}{\updefault}{\color[rgb]{0,0,0}$g_0$}%
}}}}
\put(3106,254){\makebox(0,0)[b]{\smash{{\SetFigFont{10}{12.0}{\rmdefault}{\mddefault}{\updefault}{\color[rgb]{0,0,0}$g$}%
}}}}
\put(3376,254){\makebox(0,0)[b]{\smash{{\SetFigFont{10}{12.0}{\rmdefault}{\mddefault}{\updefault}{\color[rgb]{0,0,0}$=$}%
}}}}
\put(766,254){\makebox(0,0)[b]{\smash{{\SetFigFont{10}{12.0}{\rmdefault}{\mddefault}{\updefault}{\color[rgb]{0,0,0}$=$}%
}}}}
\end{picture}%
\end{center}

\noindent One says that the $3$-polygraph $(\Sigma,R)$ \emph{terminates} if there exist no family $(f_n)_{n\in\Nb}$ of circuits of $\mon{\Sigma}$ such that $f_n\red{R}\:f_{n+1}$ for every natural number $n$.
\end{defn}

\noindent Therafter, we assume that every polygraph we consider has one $0$-cell and one $1$-cell. As for any kind of rewriting system, the easiest way to prove that a $3$-polygraph terminates is to produce a well-chosen termination order.

\begin{defn}
A \emph{termination order} on a $2$-polygraph $\Sigma$ is a strict order $>$ on parallel circuits such that there exist no family $(f_n)_{n\in\Nb}$ of circuits with $f_n>f_{n+1}$ for every $n$ and such that, for any circuit $f$, the maps $f\ast_0(\cdot)$, $(\cdot)\ast_0 f$, $f\ast_1(\cdot)$ and $(\cdot)\ast_1 f$ are strictly monotone.
\end{defn}

\begin{prop}\label{ordre_reduction}
Let $(\Sigma,R)$ be a $3$-polygraph and $>$ be a termination order on $\Sigma$. If, for any $3$-cell $\alpha$ from~$f$ to~$g$, the inequality $f>g$ holds, then $(\Sigma,R)$ terminates.
\end{prop}

\section{Crafting termination orders for $\mathbf{3}$-polygraphs}

\noindent Proposition \ref{ordre_reduction} would remain useless without a recipe to build termination orders, such as the ones that exist for term rewriting. Moreover, even though circuits are deeply linked with terms, there exist obstructions to directly transpose techniques from term rewriting to polygraphs. However, it is possible to adapt them.

Let us give the rough idea. Given a $2$-polygraph $\Sigma$, circuits of $\mon{\Sigma}$ are compared according to the "heat" they produce when presented with some "courant intensities". The courants are plugged into each input and each output of a given circuit $f$. Then, they propagate through $f$ to reach all the circuit components (elements of $\Sigma$) used to build $f$. Each component produces some heat, depending on the intensities of the courants it receives. The heat produced by $f$ is the sum of all the heats produced by the components of $f$. Given another circuit $g$, parallel to $f$, $f$ will be declared greater than $g$ if it always produces more heat than $g$ when both receive the same courant intensities. 

In order to formalize these ideas, we use two non-empty ordered sets $X$ and $Y$, for the courants:~$X$ is for descending courants, or courants going from the inputs to the outputs, and $Y$ for ascending courants. We need also a commutative monoid $M$, equipped with an order relation, such that the sum is strictly monotone in both arguments: this is used to express heats. Finally, for each $2$-cell $\phi$ in $\Sigma$, we require three \emph{monotone} maps $\phi_*:X^m\fl X^n$, $\phi^*:Y^n\fl Y^m$ and $[\phi]:X^m\times Y^n\fl M$, respectively expressing how $\phi$ transmits descending courants, how $\phi$ transmits ascending courants and how much heat it produces.

\begin{defn}\label{extension}
The three interpretations $(\cdot)_*$, $(\cdot)^*$ and $[\cdot]$ are extended from $2$-cells to circuits this way:
$$
\begin{array}{r c l c r c l c r c l}
n_* &=& \Id_{X^n} &\quad& n^* &=& \Id_{Y^n} &\quad& [n](\vec{x},\vec{y}) &=& 0 \\
(f\star_0g)_* &=& (f_*,g_*) && (f\star_0g)^* &=& (f^*,g^*) && [f\star_0g](\vec{x},\vec{x}',\vec{y},\vec{y}') &=& [f](\vec{x},\vec{y})+[g](\vec{x}',\vec{y}') \\
(f\star_1g)_* &=& g_*\circ f_* && (f\star_1g)^* &=& f^*\circ g^* &&
[f\star_1 g](\vec{x},\vec{y}) &=& [f](\vec{x},g^*(\vec{y}))+[g](f_*(\vec{x}),\vec{y})
\end{array}
$$
\end{defn}

\noindent One has to prove that the three interpretations $(\cdot)_*$, $(\cdot)^*$ and $[\cdot]$ are well-defined on every circuit and that, for each circuit $f$, the three maps $f_*$, $f^*$ and $[f]$ are monotone [Guiraud 2004]. Now we define an order on parallel circuits and prove the main result.

\begin{defn}\label{construction}
With the same notations, one defines a binary relation $>$ on parallel circuits of $\mon{\Sigma}$: let~$f$ and $g$ be two circuits with $m$ inputs and $n$ outputs. Then $f>g$ if, for any $\vec{x}\in X^m$, $\vec{y}\in Y^n$, the inequalities $f_*(\vec{x})\geq g_*(\vec{x})$, $f^*(\vec{y})\geq g^*(\vec{y})$ and $[f](\vec{x},\vec{y})>[g](\vec{x},\vec{y})$ hold.
\end{defn}

\begin{thm}\label{principal}
Let us keep the aforegiven notations and let us assume that the order relation on the commutative monoid $M$ does not admit infinite strictly decreasing sequences. Then, the binary relation~$>$ on parallel circuits of $\mon{\Sigma}$ is a termination order on $\Sigma$. In particular, if every $3$-cell $\alpha$ in $R$ from $f$ to $g$ satisfies $f>g$, then the $3$-polygraph $(\Sigma,R)$ terminates.
\end{thm}

\section{Termination orders at work}

The theorem \ref{principal} has been used in [Guiraud 2004] in order to prove two conjectures from [Lafont 2003]. We present one of them here: it states the termination of the $3$-polygraph $\lz$, which is a presentation of the structure of $\zb$-vector space. This is an important point for polygraphs since such a presentation cannot exist in the term rewriting formalism. 

\begin{figure}[!h]
\includegraphics{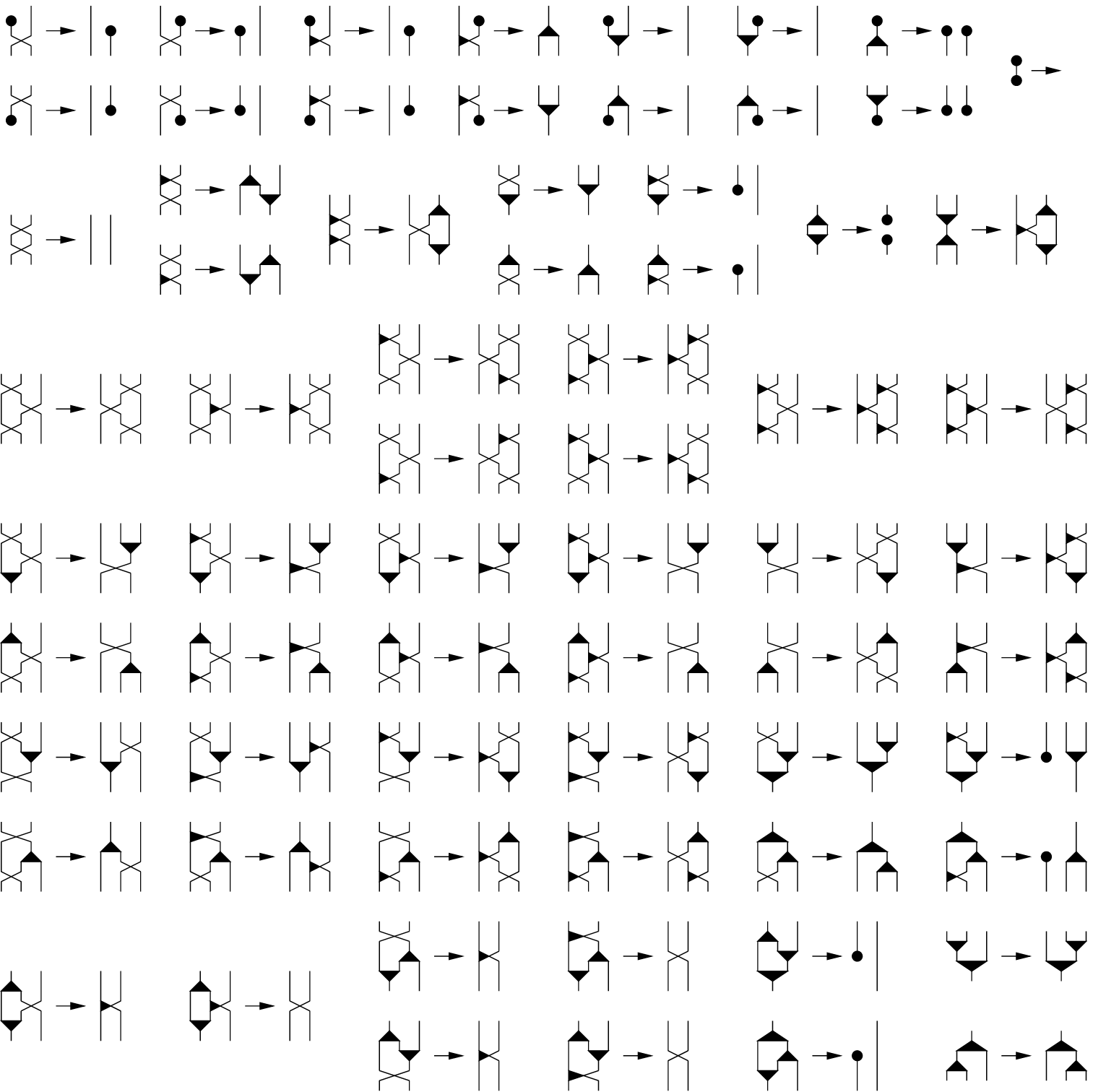}
\caption{The sixty-seven $3$-cells of $\lz$}\label{fig:3-cellules-lz2}
\end{figure}

The polygraph $\lz$ has six $2$-cells \includegraphics{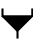},\includegraphics{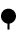}, \includegraphics{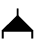}, \includegraphics{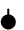}, \includegraphics{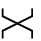} and \includegraphics{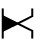}, together with sixty-seven $3$-cells, pictured in figure \ref{fig:3-cellules-lz2}. In order to prove the termination of $\lz$, we consider $X=Y=\Nb$, equipped with its natural order, while $M$ is the free commutative monoid generated by $\Nb^*$, equipped with the \emph{multiset order}: this is the smaller order strictly compatible with the sum such that $p.\ul{n}\:<\:\ul{n+1}$, for every $p$ and $n$ and where $\ul{n}$ denotes the natural number $n$ seen as a generator of $M$.

An application of theorem \ref{principal} shows that the following interpretations generate a termination order that proves the conjecture. For each $2$-cell $\alpha$, the first two diagrams give $\alpha_*$ and $\alpha^*$, while the third one gives $[\alpha]$:

\medskip
\begin{center}
\begin{picture}(0,0)%
\includegraphics{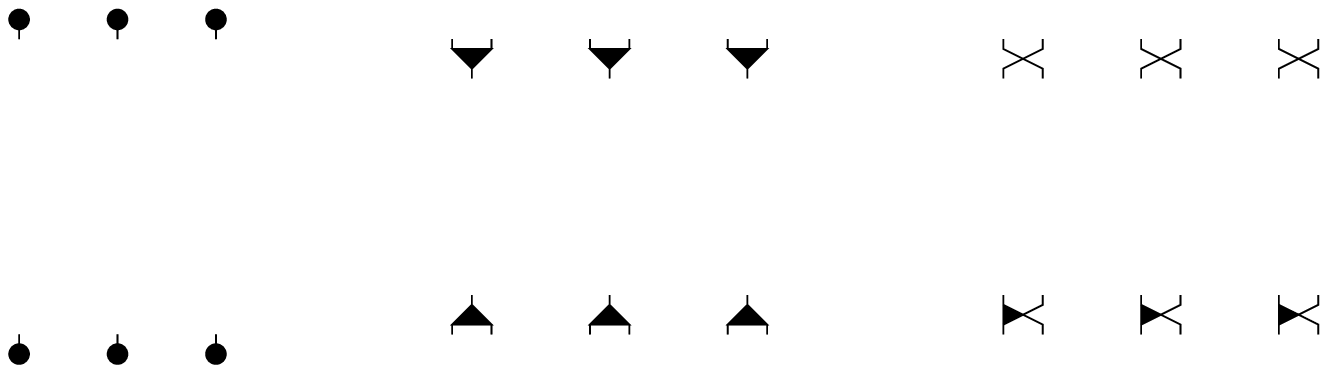}%
\end{picture}%
\setlength{\unitlength}{4144sp}%
\begingroup\makeatletter\ifx\SetFigFont\undefined%
\gdef\SetFigFont#1#2#3#4#5{%
  \reset@font\fontsize{#1}{#2pt}%
  \fontfamily{#3}\fontseries{#4}\fontshape{#5}%
  \selectfont}%
\fi\endgroup%
\begin{picture}(6325,2145)(-204,-1366)
\put(5986,-916){\makebox(0,0)[b]{\smash{{\SetFigFont{10}{12.0}{\rmdefault}{\mddefault}{\updefault}{\color[rgb]{0,0,0}$\leadsto$}%
}}}}
\put(1891,479){\makebox(0,0)[b]{\smash{{\SetFigFont{10}{12.0}{\rmdefault}{\mddefault}{\updefault}{\color[rgb]{0,0,0}$\downarrow$}%
}}}}
\put(2071,479){\makebox(0,0)[b]{\smash{{\SetFigFont{10}{12.0}{\rmdefault}{\mddefault}{\updefault}{\color[rgb]{0,0,0}$\downarrow$}%
}}}}
\put(1891,659){\makebox(0,0)[b]{\smash{{\SetFigFont{10}{12.0}{\rmdefault}{\mddefault}{\updefault}$i$}}}}
\put(2071,659){\makebox(0,0)[b]{\smash{{\SetFigFont{10}{12.0}{\rmdefault}{\mddefault}{\updefault}$j$}}}}
\put(1981, 29){\makebox(0,0)[b]{\smash{{\SetFigFont{10}{12.0}{\rmdefault}{\mddefault}{\updefault}{\color[rgb]{0,0,0}$\downarrow$}%
}}}}
\put(1981,-151){\makebox(0,0)[b]{\smash{{\SetFigFont{10}{12.0}{\rmdefault}{\mddefault}{\updefault}$i+j$}}}}
\put(3151,659){\makebox(0,0)[b]{\smash{{\SetFigFont{10}{12.0}{\rmdefault}{\mddefault}{\updefault}$i$}}}}
\put(3331,659){\makebox(0,0)[b]{\smash{{\SetFigFont{10}{12.0}{\rmdefault}{\mddefault}{\updefault}$j$}}}}
\put(3331,479){\makebox(0,0)[b]{\smash{{\SetFigFont{10}{12.0}{\rmdefault}{\mddefault}{\updefault}{\color[rgb]{0,0,0}$\downarrow$}%
}}}}
\put(3151,479){\makebox(0,0)[b]{\smash{{\SetFigFont{10}{12.0}{\rmdefault}{\mddefault}{\updefault}{\color[rgb]{0,0,0}$\downarrow$}%
}}}}
\put(2701,479){\makebox(0,0)[b]{\smash{{\SetFigFont{10}{12.0}{\rmdefault}{\mddefault}{\updefault}{\color[rgb]{0,0,0}$\uparrow$}%
}}}}
\put(2521,659){\makebox(0,0)[b]{\smash{{\SetFigFont{10}{12.0}{\rmdefault}{\mddefault}{\updefault}$i$}}}}
\put(2701,659){\makebox(0,0)[b]{\smash{{\SetFigFont{10}{12.0}{\rmdefault}{\mddefault}{\updefault}$i$}}}}
\put(2611, 29){\makebox(0,0)[b]{\smash{{\SetFigFont{10}{12.0}{\rmdefault}{\mddefault}{\updefault}{\color[rgb]{0,0,0}$\uparrow$}%
}}}}
\put(2521,479){\makebox(0,0)[b]{\smash{{\SetFigFont{10}{12.0}{\rmdefault}{\mddefault}{\updefault}{\color[rgb]{0,0,0}$\uparrow$}%
}}}}
\put(2611,-151){\makebox(0,0)[b]{\smash{{\SetFigFont{10}{12.0}{\rmdefault}{\mddefault}{\updefault}$i$}}}}
\put(3241, 29){\makebox(0,0)[b]{\smash{{\SetFigFont{10}{12.0}{\rmdefault}{\mddefault}{\updefault}{\color[rgb]{0,0,0}$\uparrow$}%
}}}}
\put(3241,-151){\makebox(0,0)[b]{\smash{{\SetFigFont{10}{12.0}{\rmdefault}{\mddefault}{\updefault}$k$}}}}
\put(-89, 29){\makebox(0,0)[b]{\smash{{\SetFigFont{10}{12.0}{\rmdefault}{\mddefault}{\updefault}$1$}}}}
\put(-89,209){\makebox(0,0)[b]{\smash{{\SetFigFont{10}{12.0}{\rmdefault}{\mddefault}{\updefault}{\color[rgb]{0,0,0}$\downarrow$}%
}}}}
\put(361,209){\makebox(0,0)[b]{\smash{{\SetFigFont{10}{12.0}{\rmdefault}{\mddefault}{\updefault}{\color[rgb]{0,0,0}$\uparrow$}%
}}}}
\put(811,209){\makebox(0,0)[b]{\smash{{\SetFigFont{10}{12.0}{\rmdefault}{\mddefault}{\updefault}{\color[rgb]{0,0,0}$\uparrow$}%
}}}}
\put(361, 29){\makebox(0,0)[b]{\smash{{\SetFigFont{10}{12.0}{\rmdefault}{\mddefault}{\updefault}$i$}}}}
\put(811, 29){\makebox(0,0)[b]{\smash{{\SetFigFont{10}{12.0}{\rmdefault}{\mddefault}{\updefault}$i$}}}}
\put(4411,479){\makebox(0,0)[b]{\smash{{\SetFigFont{10}{12.0}{\rmdefault}{\mddefault}{\updefault}{\color[rgb]{0,0,0}$\downarrow$}%
}}}}
\put(4411, 29){\makebox(0,0)[b]{\smash{{\SetFigFont{10}{12.0}{\rmdefault}{\mddefault}{\updefault}{\color[rgb]{0,0,0}$\downarrow$}%
}}}}
\put(4591,479){\makebox(0,0)[b]{\smash{{\SetFigFont{10}{12.0}{\rmdefault}{\mddefault}{\updefault}{\color[rgb]{0,0,0}$\downarrow$}%
}}}}
\put(4411,659){\makebox(0,0)[b]{\smash{{\SetFigFont{10}{12.0}{\rmdefault}{\mddefault}{\updefault}$i$}}}}
\put(4591,659){\makebox(0,0)[b]{\smash{{\SetFigFont{10}{12.0}{\rmdefault}{\mddefault}{\updefault}$j$}}}}
\put(5671,659){\makebox(0,0)[b]{\smash{{\SetFigFont{10}{12.0}{\rmdefault}{\mddefault}{\updefault}$i$}}}}
\put(5851,659){\makebox(0,0)[b]{\smash{{\SetFigFont{10}{12.0}{\rmdefault}{\mddefault}{\updefault}$j$}}}}
\put(5851,479){\makebox(0,0)[b]{\smash{{\SetFigFont{10}{12.0}{\rmdefault}{\mddefault}{\updefault}{\color[rgb]{0,0,0}$\downarrow$}%
}}}}
\put(5671,479){\makebox(0,0)[b]{\smash{{\SetFigFont{10}{12.0}{\rmdefault}{\mddefault}{\updefault}{\color[rgb]{0,0,0}$\downarrow$}%
}}}}
\put(5041,479){\makebox(0,0)[b]{\smash{{\SetFigFont{10}{12.0}{\rmdefault}{\mddefault}{\updefault}{\color[rgb]{0,0,0}$\uparrow$}%
}}}}
\put(5221,479){\makebox(0,0)[b]{\smash{{\SetFigFont{10}{12.0}{\rmdefault}{\mddefault}{\updefault}{\color[rgb]{0,0,0}$\uparrow$}%
}}}}
\put(4951,659){\makebox(0,0)[b]{\smash{{\SetFigFont{10}{12.0}{\rmdefault}{\mddefault}{\updefault}$i+j$}}}}
\put(5266,659){\makebox(0,0)[b]{\smash{{\SetFigFont{10}{12.0}{\rmdefault}{\mddefault}{\updefault}$i$}}}}
\put(4321,-151){\makebox(0,0)[b]{\smash{{\SetFigFont{10}{12.0}{\rmdefault}{\mddefault}{\updefault}$i+j$}}}}
\put(4636,-151){\makebox(0,0)[b]{\smash{{\SetFigFont{10}{12.0}{\rmdefault}{\mddefault}{\updefault}$i$}}}}
\put(4591, 29){\makebox(0,0)[b]{\smash{{\SetFigFont{10}{12.0}{\rmdefault}{\mddefault}{\updefault}{\color[rgb]{0,0,0}$\downarrow$}%
}}}}
\put(5041, 29){\makebox(0,0)[b]{\smash{{\SetFigFont{10}{12.0}{\rmdefault}{\mddefault}{\updefault}{\color[rgb]{0,0,0}$\uparrow$}%
}}}}
\put(5221, 29){\makebox(0,0)[b]{\smash{{\SetFigFont{10}{12.0}{\rmdefault}{\mddefault}{\updefault}{\color[rgb]{0,0,0}$\uparrow$}%
}}}}
\put(5041,-151){\makebox(0,0)[b]{\smash{{\SetFigFont{10}{12.0}{\rmdefault}{\mddefault}{\updefault}$i$}}}}
\put(5221,-151){\makebox(0,0)[b]{\smash{{\SetFigFont{10}{12.0}{\rmdefault}{\mddefault}{\updefault}$j$}}}}
\put(5851, 29){\makebox(0,0)[b]{\smash{{\SetFigFont{10}{12.0}{\rmdefault}{\mddefault}{\updefault}{\color[rgb]{0,0,0}$\uparrow$}%
}}}}
\put(5671, 29){\makebox(0,0)[b]{\smash{{\SetFigFont{10}{12.0}{\rmdefault}{\mddefault}{\updefault}{\color[rgb]{0,0,0}$\uparrow$}%
}}}}
\put(5671,-151){\makebox(0,0)[b]{\smash{{\SetFigFont{10}{12.0}{\rmdefault}{\mddefault}{\updefault}$k$}}}}
\put(5851,-151){\makebox(0,0)[b]{\smash{{\SetFigFont{10}{12.0}{\rmdefault}{\mddefault}{\updefault}$l$}}}}
\put(991,434){\makebox(0,0)[b]{\smash{{\SetFigFont{10}{12.0}{\rmdefault}{\mddefault}{\updefault}{\color[rgb]{0,0,0}$\leadsto$}%
}}}}
\put(1126,434){\makebox(0,0)[lb]{\smash{{\SetFigFont{10}{12.0}{\rmdefault}{\mddefault}{\updefault}$\ul{i}$}}}}
\put(3466,254){\makebox(0,0)[b]{\smash{{\SetFigFont{10}{12.0}{\rmdefault}{\mddefault}{\updefault}{\color[rgb]{0,0,0}$\leadsto$}%
}}}}
\put(3601,254){\makebox(0,0)[lb]{\smash{{\SetFigFont{10}{12.0}{\rmdefault}{\mddefault}{\updefault}$\ul{i}+\ul{k}$}}}}
\put(5986,254){\makebox(0,0)[b]{\smash{{\SetFigFont{10}{12.0}{\rmdefault}{\mddefault}{\updefault}{\color[rgb]{0,0,0}$\leadsto$}%
}}}}
\put(6121,254){\makebox(0,0)[lb]{\smash{{\SetFigFont{10}{12.0}{\rmdefault}{\mddefault}{\updefault}$\ul{i}+\ul{k}$}}}}
\put(1981,-691){\makebox(0,0)[b]{\smash{{\SetFigFont{10}{12.0}{\rmdefault}{\mddefault}{\updefault}{\color[rgb]{0,0,0}$\downarrow$}%
}}}}
\put(1891,-1141){\makebox(0,0)[b]{\smash{{\SetFigFont{10}{12.0}{\rmdefault}{\mddefault}{\updefault}{\color[rgb]{0,0,0}$\downarrow$}%
}}}}
\put(1981,-511){\makebox(0,0)[b]{\smash{{\SetFigFont{10}{12.0}{\rmdefault}{\mddefault}{\updefault}$i$}}}}
\put(1891,-1321){\makebox(0,0)[b]{\smash{{\SetFigFont{10}{12.0}{\rmdefault}{\mddefault}{\updefault}$i$}}}}
\put(2071,-1141){\makebox(0,0)[b]{\smash{{\SetFigFont{10}{12.0}{\rmdefault}{\mddefault}{\updefault}{\color[rgb]{0,0,0}$\downarrow$}%
}}}}
\put(2071,-1321){\makebox(0,0)[b]{\smash{{\SetFigFont{10}{12.0}{\rmdefault}{\mddefault}{\updefault}$i$}}}}
\put(3241,-691){\makebox(0,0)[b]{\smash{{\SetFigFont{10}{12.0}{\rmdefault}{\mddefault}{\updefault}{\color[rgb]{0,0,0}$\downarrow$}%
}}}}
\put(2611,-691){\makebox(0,0)[b]{\smash{{\SetFigFont{10}{12.0}{\rmdefault}{\mddefault}{\updefault}{\color[rgb]{0,0,0}$\uparrow$}%
}}}}
\put(2611,-511){\makebox(0,0)[b]{\smash{{\SetFigFont{10}{12.0}{\rmdefault}{\mddefault}{\updefault}$i+j$}}}}
\put(2521,-1141){\makebox(0,0)[b]{\smash{{\SetFigFont{10}{12.0}{\rmdefault}{\mddefault}{\updefault}{\color[rgb]{0,0,0}$\uparrow$}%
}}}}
\put(2701,-1141){\makebox(0,0)[b]{\smash{{\SetFigFont{10}{12.0}{\rmdefault}{\mddefault}{\updefault}{\color[rgb]{0,0,0}$\uparrow$}%
}}}}
\put(2521,-1321){\makebox(0,0)[b]{\smash{{\SetFigFont{10}{12.0}{\rmdefault}{\mddefault}{\updefault}$i$}}}}
\put(2701,-1321){\makebox(0,0)[b]{\smash{{\SetFigFont{10}{12.0}{\rmdefault}{\mddefault}{\updefault}$j$}}}}
\put(3151,-1141){\makebox(0,0)[b]{\smash{{\SetFigFont{10}{12.0}{\rmdefault}{\mddefault}{\updefault}{\color[rgb]{0,0,0}$\uparrow$}%
}}}}
\put(3331,-1141){\makebox(0,0)[b]{\smash{{\SetFigFont{10}{12.0}{\rmdefault}{\mddefault}{\updefault}{\color[rgb]{0,0,0}$\uparrow$}%
}}}}
\put(3151,-1321){\makebox(0,0)[b]{\smash{{\SetFigFont{10}{12.0}{\rmdefault}{\mddefault}{\updefault}$j$}}}}
\put(3331,-1321){\makebox(0,0)[b]{\smash{{\SetFigFont{10}{12.0}{\rmdefault}{\mddefault}{\updefault}$k$}}}}
\put(3241,-511){\makebox(0,0)[b]{\smash{{\SetFigFont{10}{12.0}{\rmdefault}{\mddefault}{\updefault}$i$}}}}
\put(-89,-871){\makebox(0,0)[b]{\smash{{\SetFigFont{10}{12.0}{\rmdefault}{\mddefault}{\updefault}{\color[rgb]{0,0,0}$\downarrow$}%
}}}}
\put(811,-871){\makebox(0,0)[b]{\smash{{\SetFigFont{10}{12.0}{\rmdefault}{\mddefault}{\updefault}{\color[rgb]{0,0,0}$\downarrow$}%
}}}}
\put(-89,-691){\makebox(0,0)[b]{\smash{{\SetFigFont{10}{12.0}{\rmdefault}{\mddefault}{\updefault}$i$}}}}
\put(811,-691){\makebox(0,0)[b]{\smash{{\SetFigFont{10}{12.0}{\rmdefault}{\mddefault}{\updefault}$i$}}}}
\put(361,-871){\makebox(0,0)[b]{\smash{{\SetFigFont{10}{12.0}{\rmdefault}{\mddefault}{\updefault}{\color[rgb]{0,0,0}$\uparrow$}%
}}}}
\put(361,-691){\makebox(0,0)[b]{\smash{{\SetFigFont{10}{12.0}{\rmdefault}{\mddefault}{\updefault}$1$}}}}
\put(4411,-691){\makebox(0,0)[b]{\smash{{\SetFigFont{10}{12.0}{\rmdefault}{\mddefault}{\updefault}{\color[rgb]{0,0,0}$\downarrow$}%
}}}}
\put(4411,-1141){\makebox(0,0)[b]{\smash{{\SetFigFont{10}{12.0}{\rmdefault}{\mddefault}{\updefault}{\color[rgb]{0,0,0}$\downarrow$}%
}}}}
\put(4591,-691){\makebox(0,0)[b]{\smash{{\SetFigFont{10}{12.0}{\rmdefault}{\mddefault}{\updefault}{\color[rgb]{0,0,0}$\downarrow$}%
}}}}
\put(5851,-511){\makebox(0,0)[b]{\smash{{\SetFigFont{10}{12.0}{\rmdefault}{\mddefault}{\updefault}$j$}}}}
\put(5671,-511){\makebox(0,0)[b]{\smash{{\SetFigFont{10}{12.0}{\rmdefault}{\mddefault}{\updefault}$i$}}}}
\put(4591,-511){\makebox(0,0)[b]{\smash{{\SetFigFont{10}{12.0}{\rmdefault}{\mddefault}{\updefault}$j$}}}}
\put(4411,-511){\makebox(0,0)[b]{\smash{{\SetFigFont{10}{12.0}{\rmdefault}{\mddefault}{\updefault}$i$}}}}
\put(5851,-691){\makebox(0,0)[b]{\smash{{\SetFigFont{10}{12.0}{\rmdefault}{\mddefault}{\updefault}{\color[rgb]{0,0,0}$\downarrow$}%
}}}}
\put(5671,-691){\makebox(0,0)[b]{\smash{{\SetFigFont{10}{12.0}{\rmdefault}{\mddefault}{\updefault}{\color[rgb]{0,0,0}$\downarrow$}%
}}}}
\put(4591,-1141){\makebox(0,0)[b]{\smash{{\SetFigFont{10}{12.0}{\rmdefault}{\mddefault}{\updefault}{\color[rgb]{0,0,0}$\downarrow$}%
}}}}
\put(4321,-1321){\makebox(0,0)[b]{\smash{{\SetFigFont{10}{12.0}{\rmdefault}{\mddefault}{\updefault}$i+j$}}}}
\put(4636,-1321){\makebox(0,0)[b]{\smash{{\SetFigFont{10}{12.0}{\rmdefault}{\mddefault}{\updefault}$i$}}}}
\put(5041,-691){\makebox(0,0)[b]{\smash{{\SetFigFont{10}{12.0}{\rmdefault}{\mddefault}{\updefault}{\color[rgb]{0,0,0}$\uparrow$}%
}}}}
\put(5221,-691){\makebox(0,0)[b]{\smash{{\SetFigFont{10}{12.0}{\rmdefault}{\mddefault}{\updefault}{\color[rgb]{0,0,0}$\uparrow$}%
}}}}
\put(4951,-511){\makebox(0,0)[b]{\smash{{\SetFigFont{10}{12.0}{\rmdefault}{\mddefault}{\updefault}$i+j$}}}}
\put(5266,-511){\makebox(0,0)[b]{\smash{{\SetFigFont{10}{12.0}{\rmdefault}{\mddefault}{\updefault}$i$}}}}
\put(5041,-1141){\makebox(0,0)[b]{\smash{{\SetFigFont{10}{12.0}{\rmdefault}{\mddefault}{\updefault}{\color[rgb]{0,0,0}$\uparrow$}%
}}}}
\put(5221,-1141){\makebox(0,0)[b]{\smash{{\SetFigFont{10}{12.0}{\rmdefault}{\mddefault}{\updefault}{\color[rgb]{0,0,0}$\uparrow$}%
}}}}
\put(5041,-1321){\makebox(0,0)[b]{\smash{{\SetFigFont{10}{12.0}{\rmdefault}{\mddefault}{\updefault}$i$}}}}
\put(5221,-1321){\makebox(0,0)[b]{\smash{{\SetFigFont{10}{12.0}{\rmdefault}{\mddefault}{\updefault}$j$}}}}
\put(5671,-1141){\makebox(0,0)[b]{\smash{{\SetFigFont{10}{12.0}{\rmdefault}{\mddefault}{\updefault}{\color[rgb]{0,0,0}$\uparrow$}%
}}}}
\put(5851,-1141){\makebox(0,0)[b]{\smash{{\SetFigFont{10}{12.0}{\rmdefault}{\mddefault}{\updefault}{\color[rgb]{0,0,0}$\uparrow$}%
}}}}
\put(5671,-1321){\makebox(0,0)[b]{\smash{{\SetFigFont{10}{12.0}{\rmdefault}{\mddefault}{\updefault}$k$}}}}
\put(5851,-1321){\makebox(0,0)[b]{\smash{{\SetFigFont{10}{12.0}{\rmdefault}{\mddefault}{\updefault}$l$}}}}
\put(1126,-1096){\makebox(0,0)[lb]{\smash{{\SetFigFont{10}{12.0}{\rmdefault}{\mddefault}{\updefault}$\ul{i}$}}}}
\put(991,-1096){\makebox(0,0)[b]{\smash{{\SetFigFont{10}{12.0}{\rmdefault}{\mddefault}{\updefault}{\color[rgb]{0,0,0}$\leadsto$}%
}}}}
\put(3601,-916){\makebox(0,0)[lb]{\smash{{\SetFigFont{10}{12.0}{\rmdefault}{\mddefault}{\updefault}$\ul{i}+\ul{j}$}}}}
\put(3466,-916){\makebox(0,0)[b]{\smash{{\SetFigFont{10}{12.0}{\rmdefault}{\mddefault}{\updefault}{\color[rgb]{0,0,0}$\leadsto$}%
}}}}
\put(6121,-916){\makebox(0,0)[lb]{\smash{{\SetFigFont{10}{12.0}{\rmdefault}{\mddefault}{\updefault}$\ul{i}+\ul{k}$}}}}
\end{picture}%
\end{center}

\vfill\pagebreak

\section*{References}
\emptysectionmark{References}

\noindent \textsc{A. Burroni}, \emph{Higher-dimensional word problems with applications to equational logic}, Theoretical computer science 115, 1993.

\smallskip
\noindent  \textsc{Y. Guiraud}, \emph{Termination orders for $3$-dimensional rewriting}, J. of pure and applied algebra, to appear (2004).

\smallskip
\noindent \textsc{Y. Lafont}, \emph{Towards an algebraic theory of boolean circuits}, J. of pure and applied algebra 184, 2003.

\end{document}